\documentclass{article}
\usepackage[leqno]{amsmath}
\usepackage{amsthm,amscd}
\usepackage{amsfonts} % to get the \mathbb alphabet
\usepackage{amssymb} % for subsetneq
\newtheorem{theorem}{Theorem}[section]
\newtheorem{definition}[theorem]{Definition}
\newtheorem{example}[theorem]{Example}
\newtheorem{lemma}[theorem]{Lemma}
\newtheorem{proposition}[theorem]{Proposition}
\newtheorem{corollary}[theorem]{Corollary}
\newtheorem{remark}[theorem]{Remark}

\newcommand{\bb}[1]{ {\mathbb{#1}} }
\newcommand{\A}{{\mathcal A}}

\newcommand{\bK}{{\mathbf K}}

\newcommand{\Poin}{{\rm Poin}}

\newcommand{\codim}{{\rm codim}}

\newcommand{\Der}{{\rm Der}}

\begin{document}
\title{The Poincar\'e series of the algebra of\\
rational functions which are regular outside\\ 
hyperplanes}
\author{{\sc Hiroki Horiuchi} and {\sc Hiroaki Terao}
\footnote{partially supported by 
the Grant-in-aid for scientific research (No.1144002), 
the Ministry of Education, Sports, Science and Technology, Japan
}\\
{\small \it Tokyo Metropolitan University, Mathematics Department}\\
{\small \it Minami-Ohsawa, Hachioji, Tokyo 192-0397, Japan}
}
\date{}
\maketitle

\begin{abstract}

Let $\Delta$ be a finite set of nonzero linear forms in several
variables with coefficients in a field $\mathbf K$ of characteristic
zero.
Consider the $\mathbf K$-algebra 
$R(\Delta)$ of rational functions 
on $V$ which are regular outside
$\bigcup_{\alpha\in\Delta} \ker\alpha$.
Then the ring $R(\Delta)$ is naturally doubly filtered
by the degrees of denominators and of numerators.
In this paper we give an explicit
combinatorial formula for the Poincar\'e series in two variables 
of the associated bigraded vector space $\overline{R}(\Delta) $.
This generalizes the main theorem of \cite{ter3}.

\noindent
{\it Mathematics Subject Classification (2000): 32S22, 13D40, 52C35} 
\end{abstract}

\bigskip
\setcounter{section}{0}
\setcounter{equation}{0}

\section{Introduction and main results}

\bigskip
Let $V$ be a
vector space of dimension $\ell$ over
a field $\bK$ of characteristic zero. Let $V^*$ be the dual
space of $V$. Suppose that 
$\Delta$ is a finite subset of nonzero elements
of $V^*$. We assume that no two vectors of $\Delta$ are
proportional. Let $S=S(V^*)$ be the symmetric algebra of $V^*$.
Then $S$  may be regarded as the ring of polynomial functions on
$V$.  

\begin{definition}
Let $R(\Delta)$ be the $\bK$-algebra 
 of rational functions on $V$ which are regular
outside the set
$\bigcup_{\alpha \in \Delta} \{\alpha=0\}$. 
In other words, $R(\Delta) = S\left[\Delta^{-1} \right]$,
where $\Delta^{-1} := \{1/\alpha \mid \alpha\in\Delta\}$.
\end{definition}

Let $f/g$ be a homogeneous element of $R(\Delta)$, i.e., 
both $f$ and $g$ are homogeneous polynomials. 
Then the (total) {\bf
degree} of $f/g$ is defined by $\deg(f/g) = \deg f - \deg g$. 
Let $R = R(\Delta), p, q\in {\mathbb Z}$.
Define
a {\bf K}-vector subspace
\[
R^p_q := \left\langle \frac{f}{g} \in R \mid
 \deg f \le p,\ 
 \deg g \le q
\right\rangle_{{\bf K}} 
\]
 of $R$.
 Here the notation $\left\langle  
 T\right\rangle_{\bK} $ stands for the 
 $\bK$-vector subspace of $R$  spanned by a subset $T$  of $R$.    
Agree that $R^{p}_{q}=0  $ if either $p<0$ or $q<0$. 
Then we have a double filtration on $R$:
$$\dots\supseteq 
R^p_q \supseteq  R^{p-1}_q \supseteq \dots 
\hspace{5mm} 
{\rm and}
\hspace{5mm} 
\dots\supseteq 
R^p_q \supseteq  R^{p}_{q-1} \supseteq \dots.$$
Define
\[
\overline{R^p_q} := R^p_q \left/ ( R^{p-1}_q + R^p_{q-1} ) \right. 
\ \ 
{\rm and}
\ \ 
\overline{R} = \overline{R}(\Delta) := \bigoplus_{p, q}  \overline{R^p_q}.
\]
The vector space $\overline{R^p_q}$ can be
considered, roughly speaking, as
the space of the homogeneous rational functions expressible 
as sums of fractions with numerators of degree $p$ 
and denominators of degree $q$ where $p$ and $q$ are 
the smallest possible.
Since each $\overline{R^{p}_{q}  } $
is finite-dimensional, one can define 
the {\bf Poincar\'e series}
\[
\Poin(\overline{R}(\Delta),s,t) := \sum_{p, q} 
\dim(\overline{R^p_q}) s^p t^q
\]
of the bigraded vector space
$\overline{R}(\Delta)$. 
Let
$$
\mathcal{A}(\Delta)=\{\ker(\alpha)\mid\alpha \in \Delta\}.
$$
Then $\mathcal{A}(\Delta)$ is a (central) arrangement of hyperplanes
\cite{ort1} in $V$. 
Let $\Poin(\A(\Delta), \ \cdot)$ be the Poincar\'e
polynomial 
 \cite[Definition 2.48]{ort1} 
 of the arrangement $\A(\Delta)$,
 which is combinatorially defined. 
 (See Definition \ref{poincare}.)
Our main theorem is the following
\begin{theorem}
\label{maintheorem}
We have
\[
\Poin(\overline{R}(\Delta),s,t) =  \frac{1}{(1-s)^{\ell}}
\Poin\left(\A(\Delta), \
\frac{t(1-s)}{(1-t)} \right).
\]
\end{theorem}

By combining Theorem \ref{maintheorem} and a later result
(Proposition \ref{2.2}), one has
\begin{corollary}
\label{1.2B}
\[
\sum_{p, q} \dim(R^{p}_{q}) s^{p}t^{q} =
(1-s)^{-\ell-1} (1-t)^{-1} 
\Poin\left(\A(\Delta), \
\frac{t(1-s)}{(1-t)} \right).
\]
\end{corollary}

By Theorem~\ref{maintheorem} and the factorization theorem
(Theorem \ref{factorizationtheorem})
in \cite{ter1}, we may easily show the following

\begin{corollary}
\label{4.2} 
If $\mathcal{A}(\Delta)$ is a free arrangement with exponents
$(d_1, \cdots, d_\ell)$ 
\cite{ter1} \cite[Definitions 4.15, 4.25]{ort1},
then
$$
{\rm Poin}(\overline{R}(\Delta), t)
=
(1-s)^{-\ell} 
(1-t)^{-\ell} 
\prod^{\ell}_{i=1}
\{1+(d_{i}-1-s d_{i})t\}.
$$
\end{corollary}

In particular, one has

\begin{corollary}
\label{reflectionarrangement}
When  $\A(\Delta)$ is the set of reflecting hyperplanes of 
any (real or complex) reflection group with exponents
$(d_1, \cdots, d_\ell)$, 
\cite[Definitions 6.22]{ort1},
then
$$
{\rm Poin}(\overline{R}(\Delta), t)
=
(1-s)^{-\ell} 
(1-t)^{-\ell} 
\prod^{\ell}_{i=1}
\{1+(d_{i}-1-s d_{i})t\}.
$$
\end{corollary}

\begin{example}
\label{1.6} 
{\rm
Let $x_{1},\dots, x_{\ell}   $ be a basis for $V^{*} $.
Let $\Delta = \{x_{i} - x_{j} \mid  1\leq i < j \leq \ell\}$.
Then $\Delta$ gives the root system of type $A_{\ell-1}$.
The corresponding
reflection arrangement $\A(\Delta)$ is a braid arrangement 
with exponents $(0, 1, \dots , \ell-1)$ \cite[Example 4.32]{ort1}.
So, by Corollary~\ref{reflectionarrangement}, we have
\begin{multline*}
{\rm Poin}(\overline{R}(\Delta), s, t) 
=
(1-s)^{-\ell} (1-t)^{-\ell + 1}\prod_{i=1}^{\ell-1}  
\{1 + (i-1)t - ist\}
\\
=
(1-s)^{-\ell} (1-t)^{-\ell + 1} (1-st)(1+t-2st)
\dots (1+(\ell - 2)t-(\ell-1)st).  
\end{multline*} 
For instance, when $\ell=3$, we have 
\[
R(\Delta)={\bf K}\left[x_{1}, x_{2}, x_{3},   
\frac{1}{x_{1} - x_{2}},
\frac{1}{x_{2} - x_{3}},
\frac{1}{x_{1} - x_{3}}
\right]
\]
and
 \begin{eqnarray*} 
&~&{\rm Poin}(\overline{R}(\Delta), s, t) \\
&=&
\frac{(1-st)(1+t-2st)}{(1-s)^{3} (1-t)^{2}}\\
&=&
\sum_{p \ge 0} \binom{p+2}{2} s^{p} 
+
\sum_{p \ge 0} \sum_{q > 0} (3p+2q+1) s^{p} t^{q}\\ 
&=&
1 + (3s+3t) + (6s^{2}+6st+ 5 t^{2}) + 
(10 s^{3} + 9 s^{2} t + 8 s t^{2} +7 t^{3}) \\
~~~&~& + (15 s^{4} + 12 s^{3}t + 11 s^{2} t^{2} + 10 s t^{3} + 9 t^{4}) 
+ \dots . 
\end{eqnarray*} 
Since 
the coefficient of $s t^{2} $ is equal to 8, 
the space $\overline{R^{1}_{2}}$ is eight-dimensional.
Actually we can easily verify by direct computation that
the classes of
\begin{multline*}
\frac{x_{1} }{(x_{1} - x_{2} )(x_{1} - x_{3} )}, \
\frac{x_{1} }{(x_{1} - x_{2} )(x_{2} - x_{3} )}, \
\frac{x_{1} }{(x_{1} - x_{2} )^{2}}, \
\frac{x_{3} }{(x_{1} - x_{2} )^{2}}, \ \\
\frac{x_{1} }{(x_{2} - x_{3} )^{2}}, \
\frac{x_{2} }{(x_{2} - x_{3} )^{2}}, \
\frac{x_{1} }{(x_{1} - x_{3} )^{2}}, \
\frac{x_{2} }{(x_{1} - x_{3} )^{2}}
\end{multline*}
form a basis for $\overline{R^{1}_{2}}$.
This implies that an arbitrary element $\varphi$
 of $R^{1}_{2}  $
can be uniquely expressed as a $\bK$-linear combination
of these eight rational functions modulo 
$R^{0}_{2} + R^{1}_{1}$.  
For example, when
\[
\varphi =
\frac{x_{1} + 2 x_{2} - x_{3} + 3}{(x_{1} - x_{3} )(x_{2} - x_{3} )}  
\]
$\varphi$ is decomposed into partial fractions as
\begin{eqnarray*}
\label{}
\varphi &=& 
-2\, \frac{x_{1}}{(x_{1} - x_{2} )(x_{1} - x_{3} )} 
+2\, \frac{x_{1}}{(x_{1} - x_{2} )(x_{2} - x_{3} )} 
+\frac{3}{(x_{1} - x_{3} )(x_{2} - x_{3} )} \\
&~& 
-\frac{1}{x_{2} - x_{3} }
+\frac{2}{x_{1} - x_{3} },
\end{eqnarray*}
which implies
\[
\varphi 
\equiv
-2\, \frac{x_{1}}{(x_{1} - x_{2} )(x_{1} - x_{3} )} 
+2\, \frac{x_{1}}{(x_{1} - x_{2} )(x_{2} - x_{3} )} 
\ \ 
\bmod
\ \ R^{0}_{2} + R^{1}_{1}.
\]
}
\end{example}

\begin{remark}
{\rm The specializations of
$\Poin(\overline{R}(\Delta),s,t)$ at
$s=0$ and $t=0$ have been known: 
 
(1) Suppose $t=0$. Since $\overline{R^{p}_{0}} =
R^{p}_{0}/R^{p-1}_{0}    $ is isomorphic to
the vector space of homogeneous polynomials of degree $p$ 
in $\ell$ variables,  we have
$\bigoplus_{p} \overline{R^{p}_{0}} \simeq S$ and thus
$$\Poin(\overline{R}(\Delta),s,0)=\frac{1}{(1-s)^{\ell}}.$$

(2) Suppose $s=0$.  Note that  $\overline{R^{0}_{q}} =
R^{0}_{q}/R^{0}_{q-1}$ is isomorphic to a ${\mathbf K}$-vector subspace
\[
C(\Delta)_{q} :=
\left\langle \frac{1}{g} \in R \mid g {\rm ~is~a~product~of~ } q
{\rm ~linear~forms~in~} \Delta
\right\rangle_{\bK}  
\]
 of $R(\Delta)$.  The $\bK$-subalgebra
 \[
 C(\Delta) := \bigoplus_{q} C(\Delta)_{q}
 \]
of $R(\Delta)$ was studied and the formula
$$\Poin(C(\Delta),t)
=
\sum_{q\ge 0} (\dim C(\Delta)_{q}) t^{q} 
=\Poin(\A(\Delta),  (1-t)^{-1} t)$$
was obtained in \cite[Theorem 1.4]{ter1}.
(The proof of this formula in \cite{ter1} heavily depends
on \cite{brv1} and \cite{ort2}.) 
  So we have 
$$\Poin(\overline{R}(\Delta),0,t)=\Poin(\A(\Delta), (1-t)^{-1} t).$$
}
\end{remark}

\section{The bigraded vector space $\overline{R}(\Delta)$}
Let $p, q \in {\mathbb Z}$. 
Let $\pi^p_q : R^p_q \rightarrow \overline{R^p_q}$ be
the canonical projection.
Fix a linear section  
$$\iota^p_q  : \overline{R^p_q} \longrightarrow R^p_q$$ 
for each canonical projection
$\pi^p_q$ such that
\[
{\rm Image}(\iota^p_q) \subseteq \left\langle \frac{f}{g} \in R \mid
f {\rm ~and~} g {\rm ~are~homogeneous~with~}  
 \deg f = p \,\,
{\rm and} \,\, \deg g = q
\right\rangle_{\bK}.
\]

\begin{proposition}
\label{2.1}
There exists a linear isomorphism
\[
\kappa := \sum_{p, q} \iota^p_q :
\overline{R} =
\bigoplus_{p, q} \overline{R^p_q}
 \stackrel{\sim}{\longrightarrow}
\sum_{p, q} R^p_q = R .
\]
\end{proposition}

\noindent
{\it Proof.}\/
At first, we prove the injectivity of $\kappa$.
For every $\phi  = \sum_{p, q} \phi^p_q
\in \ker(\kappa)$, 
where $\phi^p_q \in \overline{R^p_q}$,
one has
\[
0 = \kappa(\phi) =  
\sum_{p, q} \iota^p_q (\phi^p_q) =
\sum_{a \in \bb{Z}} \sum_{ p-q = a } \iota^p_q (\phi^p_q) .
\]
This implies that
$
\sum_{ p-q = a } \iota^p_q (\phi^p_q) = 0
$
for each $a \in \bb{Z}$
because $\deg \iota^p_q (\phi^p_q) = p-q$.
Suppose $\phi\neq 0$. 
Then there exists an integer $b\in \bb{Z}$ such that
$
\sum_{ p-q = b } \phi^p_q \neq 0.
$
Let $p_b := \max\{ p \mid \phi^p_q \neq 0, \,\,\,p-q = b \}$ and
$q_b = p_b - b.$
Then $\psi^{p_{b}}_{q_{b}} := \iota^{p_{b}}_{q_{b}} 
(\phi^{p_{b}}_{q_{b}}) \neq 0$,
because $\iota^{p_{b}}_{q_{b}}$ is injective.
Therefore one has
\[
0
=
\sum_{ p-q = b } \iota^p_q (\phi^p_q) = 
\psi^{p_b}_{q_b} + \psi^{p_b -1}_{q_b -1} + \cdots 
\]
and
\[
 \psi^{p_b}_{q_b} = - \psi^{p_b -1}_{q_b -1} - 
\psi^{p_b -1}_{q_b -1}
-
\cdots \ \in R^{p_b -1}_{q_b -1}
\subseteq
R^{p_b}_{q_b -1}
+
R^{p_b -1}_{q_b}.
\]
So we obtain 
 $$
 \phi^{p_b}_{q_b}
 =
\pi^{p_b}_{q_b}( \iota^{p_b}_{q_b} (\phi^{p_b}_{q_b}) ) 
= \pi^{p_b}_{q_b}(\psi^{p_b}_{q_b}) =
0,$$
which contradicts the definition of $p_b$.
Therefore we have $\ker(\kappa) = \{ \mathbf{0} \}$.

Next we prove the surjectivity of $\kappa$.
Let $p\ge 0$ and $q\ge 0$. 
Define
\[
\kappa^{p}_{q} := \sum_{a\le p, b\le q} \iota^a_b :
\bigoplus_{a\le p, b\le q} \overline{R^a_b}
{\longrightarrow}
R^p_q.
\] 
It is sufficient to show that
$\kappa^{p}_{q}$ is surjective for
$p, q \in {\mathbb Z}$.
If either $p$ or $q$ is negative, then
$R_{q}^{p} = 0$ and there is nothing to prove.
We will prove 
$R_{q}^{p} = {\rm Image}(\kappa^{p}_{q})$ for 
$p + q \ge 0$ by an induction on $p+q$.
When $p=q=0$,
$\pi_{0}^{0} : \bK = R^0_0  \rightarrow  \overline{R^0_0} = \bK$
is an isomorphism.  So $ R^{0}_{0} 
= {\rm Image}(\iota_{0}^{0}) 
= {\rm Image}(\kappa_{0}^{0})$.   
Suppose $p+q>0$.
For $
\psi \in R^p_q$,
let $\eta := \psi - \iota^{p}_{q}\circ \pi^{p}_{q}(\psi)$.
Then 
\[
\pi^{p}_{q}  (\eta) = \pi^{p}_{q}  (\psi)
 - \pi^{p}_{q}  \circ \iota^{p}_{q}\circ \pi^{p}_{q}(\psi)=0.
\]
So $\eta\in 
\ker\pi^{p}_{q} = 
R^{p-1}_{q} + R^{p}_{q-1}\subseteq {\rm Image}(\kappa^{p}_{q}) $
by the induction assumption.
We thus have
\[
\psi = \eta + \iota^{p}_{q}\circ \pi^{p}_{q}(\psi) \in 
 {\rm Image}(\kappa^{p}_{q}),
 \]
which implies that $\kappa^{p}_{q}$ is surjective. ~~~~$\square$
\bigskip

In the course of the proof of Proposition \ref{2.1}, 
we have already proved the following  

\begin{proposition}
\label{2.2} 
Let $p\ge 0$ and $q\ge 0$.  Then
there exists a linear isomorphism
\[
\kappa^{p}_{q} = \sum_{a\le p, b\le q} \iota^a_b :
\bigoplus_{a\le p, b\le q} \overline{R^a_b}
 \stackrel{\sim}{\longrightarrow}
R^p_q.
\] 
\end{proposition}

\section{A decomposition of $\overline{R^{p}_{q}  } $ }

Let ${\bf E}_{p}(\Delta)$ be the set of all $p$-tuples composed of
elements of $\Delta$. Let    
${\bf E}(\Delta) := \bigcup_{p\geq 0} {\bf E}_{p}(\Delta)$.
The union is disjoint.
Write $\prod {\cal E} := 
\alpha_{1} \dots \alpha_{p} \in S$ when ${\mathcal E}=
(\alpha_{1}, \dots, \alpha_{p})\in {\bf E}_{p}(\Delta)   $.
For $\varepsilon \in \mathbf{E}(\Delta)$, let $V(\varepsilon)$
denote the set of common zeros of $\varepsilon$ : 
$$V(\varepsilon)=\bigcap^p_{i=1} \ker (\alpha_i)$$ 
when
$\varepsilon=(\alpha_1, \cdots, \alpha_p)$. 
Define
$$
L=L(\Delta)=\{V(\varepsilon) \mid \varepsilon \in \mathbf{E}
(\Delta)\}.
$$
Agree that $V(\varepsilon)=V$ if $\varepsilon$ is the empty tuple.
Introduce a partial order $\le$ into $L$ by reverse inclusion:
$X\le Y \Leftrightarrow X\supseteq Y$. 
Then $L$ is equal to the intersection lattice of the arrangement
$\mathcal{A}(\Delta)$ \cite[Definition 2.1]{ort1}. 
For $X\in L$, define
$$
\mathbf{E}_X(\Delta):=\{\varepsilon \in \mathbf{E}(\Delta)
\mid V(\varepsilon)=X \}. 
$$
Then
$$
\mathbf{E}(\Delta)=\bigcup_{X \in L}\mathbf{E}_X(\Delta)
~~~\text{(disjoint)}.
$$
Define
$$
C_X(\Delta):=\sum_{\varepsilon \in \mathbf{E}_X(\Delta)} \mathbf{K}
{(\prod \varepsilon)}^{-1}.
$$
The following proposition is in \cite[Proposition 2.1]{ter1}:

\begin{proposition}
\label{3.1}
$$
C(\Delta)=\bigoplus_{X\in L} C_X(\Delta).
$$
\end{proposition}

Let $X\in L$ and $q\in {\mathbb Z}_{\ge 0} $.
Define 
\[
C_{q, X}:=
C_{q, X}(\Delta) := C(\Delta)_{q} \cap C_{X}(\Delta)
=  
\sum_{\varepsilon \in \mathbf{E}_q(\Delta)\cap  \mathbf{E}_X(\Delta)} 
\mathbf{K}
{(\prod \varepsilon)}^{-1}.
\]
 Then 
 \[
 C_{X}(\Delta)
 =
 \bigoplus_{q\ge 0} C_{q, X}.
 \]
 \begin{definition}
 (\cite[2.42]{ort1})
 \label{moebius} 
 Let the 
 {\bf M\"obius function}
 \[
 \mu : L(\Delta) \rightarrow {\mathbb Z}
 \]
 be characterized by $\mu(V)=1$ 
 and for $X < V$ by 
 $
 \sum_{Y\le X} \mu(Y) = 0.
 $
 \end{definition}
 The following result is in \cite[Theorem 1.4]{ter1}:
 \begin{proposition}
\label{3.1A} 
$$
\Poin(C_{X}(\Delta), t) :=
\sum_{q} \left(\dim C_{q, X} \right) t^{q} 
=
(-1)^{\codim X} \mu(X) \left(\frac{t}{1-t}\right)^{\codim X}.
$$
\end{proposition}
 
Define $I(X) := \{ f \in S \mid f|_X \equiv 0 \}$,
i.e., $I(X)$ is the prime ideal of $S$ 
generated by the polynomial functions vanishing on $X$.
Let $S_X := S / I(X)$.
Let $S^p$ denote the degree $p$ homogeneous part of $S$.
Then $S = \bigoplus_{p \ge 0} S^p$.
Define
\[
S^p_X := S^p / S^{p} \cap I(X).
\]
Then
\[
S_X := \bigoplus_{p\ge 0} S^p_X.
\]

\begin{lemma}
\label{3.2} 
For $p\ge 0$, $q\ge 0$, and $X \in L$, 
suppose that 
$h_{1}, h_{2}, \dots , h_{m}\in C_{q, X}$ 
are linearly independent over $\bK$ and 
$f_{1}, f_{2}, \dots , f_{m}\in S^{p}$.
Then the following three conditions are
equivalent:

(1) Each 
$f_{i}$ belongs to the ideal $I(X)$,

(2) 
$
\sum_{i=1}^{m} f_{i}\, h_{i} \in R^{p-1}_{q-1},
$

(3) 
$
\sum_{i=1}^{m} f_{i} \, h_{i} \in 
R^{p-1}_{q-1}
+\sum_{{Y\in L}\atop{Y\neq X}} S^{p}\, C_{q, Y}.
$
\end{lemma}

\noindent
{\it Proof.}
(1) $\Rightarrow$ (2) :
If $f_{i} \in S^{p} \cap I(X)$, then
one can easily see that
$f_{i} \, h_{i} \in R^{p-1}_{q-1}$ for each $i$.  

(2) $\Rightarrow$ (3) : Obvious.

(3) $\Rightarrow$ (1) : 
Let 
$$
\prod\Delta = \prod_{\alpha\in\Delta} \alpha,
\ \
\Delta_{X} = \Delta \cap I(X),
\ \
{\rm and}
\ \
\prod\Delta_{X}  = \prod_{\alpha\in\Delta_{X} } \alpha.
$$ 
 Suppose that $\dim X = k$. 
Choose a basis $x_{1}, x_{2}, \dots, x_{\ell}   $ for $V^{*} $ such that
$X=V(x_{k+1}, \dots , x_{\ell})$. 
Let
\[
h'_i :=  \left( \prod \Delta_X \right)^q  h_i  
\ \in \bK[x_{k+1}, \cdots, x_{\ell}],
\]
which is a homogeneous polynomial 
of degree $q | \Delta_X | - q$.
Multiply $(\prod\Delta_{X})^{q} $ to
$
\sum_{i=1}^{m} f_{i} \, h_{i} \in R^{p-1}_{q-1}
+\sum_{{Y\in L}\atop{Y\neq X}} S^{p}\, C_{q, Y}
$
to get
\[
\left( \frac{ \prod \Delta }{ \prod \Delta_X } \right)^q \sum_i f_i 
\, h'_i 
=
\left(\prod\Delta\right)^{q}
\sum_i f_i \,
h_i 
\in {I(X)}^{q |\Delta_X| - q +1}.
\]
Since ${I(X)}^{q |\Delta_X| - q +1}$ is primary and 
$\prod \Delta / \prod \Delta_X \not\in I(X)$, one has
\[
\sum_i f_i \, h'_i 
\in {I(X)}^{q |\Delta_X| - q +1}
= (x_{k+1}, \cdots, x_{\ell})^{q |\Delta_X| - q +1}.
\]
Define
$f'_{i} := f_{i} (x_{1}, \dots, x_{k}, 0, \dots ,0)
\in \bK[x_{1}, \dots, x_{k}]$. 
Since $f_{i} - f'_{i} \in I(X)$,
one has
\[
\sum_i f'_i \, h'_i 
\in {I(X)}^{q |\Delta_X| - q +1}
= (x_{k+1}, \cdots, x_{\ell})^{q |\Delta_X| - q +1}.
\]
Recalling that each $f'_{i}$ lies in $\bK[x_{1}, \cdots, x_{k}]$
and $\deg h'_{i} = q |\Delta_X| - q$
with $h'_{i} \in 
(x_{k+1}, \cdots, x_{\ell})^{q |\Delta_X| - q}$, 
we may conclude
\[
\sum_i f'_i \, h'_i = 0.
\]
We may deduce that $f'_{1}=f'_{2} =\dots =f'_{m} = 0 $ 
because $
h'_{1},
\dots
h'_{m}
\in \bK[x_{k+1}, \cdots, x_{\ell}]
$ 
are linearly independent over $\bK$ and thus over $\bK[x_{1}, \cdots, x_{k}]$ 
also.
Therefore $f_{i}\in I(X)$ for each $i$.
$\square$

\begin{proposition}
\label{conj_main}
For $p\ge 0$, $q\ge 0$, and $X \in L$, let
\[
\tau^p_{q, X} : S^p_X \otimes C_{q, X} \longrightarrow \overline{R^p_q}
\]
be the linear map characterized by
\[
\tau^{p}_{q, X} \left([f] \otimes (\prod \varepsilon)^{-1}\right) = 
\left[\frac{f}{(\prod \varepsilon)} \right],
\]
where $f\in S^{p}$ and 
$\varepsilon \in \mathbf{E}_q(\Delta)\cap  
\mathbf{E}_X(\Delta)$. 
Then 
$\tau^{p}_{q, X}$ is injective for each $X\in L$.
\end{proposition}

\noindent
{\it Proof.} Let
$\varepsilon \in \mathbf{E}_q(\Delta)\cap  
\mathbf{E}_X(\Delta)$.  Note that $\varepsilon$
generates the ideal $I(X)$.
If $f$ belongs to 
$S^{p} \cap I(X)$, then  
$$\frac{f}{(\prod \varepsilon)}\in R^{p-1}_{q-1}\subset
R^{p-1}_{q}
+
R^{p}_{q-1}$$
by Lemma~\ref{3.2}. 
This implies that the map $\tau^{p}_{q, X}$ is well-defined.

Choose  a $\bK$-basis 
$
h_{1},
h_{2},
\dots
h_{m}$ 
for $C_{q, X}$.
Then an arbitrary element $\varphi\in S^{p}_{X} \otimes C_{q, X}  $ 
can be expressed as
$$\varphi = \sum_{i=1}^{m} \left[f_{i} \right] \otimes h_{i}  $$
for some $f_{i} \in S^{p} 
 \,\,\, (i=1,\dots,m)$.  
Suppose $\varphi\in\ker(\tau^{p}_{q, X})$,
i.e., $\sum_{i=1}^{m} f_{i}\, h_{i} \in R^{p-1}_{q} + R^{p}_{q-1}$.
This implies    $\sum_{i=1}^{m}  f_{i}\, h_{i} 
\in R^{p-1}_{q-1}.$
By Lemma~\ref{3.2}, one has $f_{i} \in I(X)$ for each
$i$.  
Therefore $\varphi=0$
and thus $\tau^{p}_{q, X}$ is injective.
$\square$

\bigskip

\begin{proposition}
\label{3.4} 
Define
$\overline{R^p_{q, X}} :=
{\rm Image}(\tau^{p}_{q, X}).$ 
Then we have
\[
\overline{R^p_q} = \bigoplus_{X \in L} \overline{R^p_{q, X}}.
\]
\end{proposition}

\noindent
{\it Proof.}
An arbitrary element $\varphi \in \overline{R^p_q}$ 
can be expressed as
$$\varphi = 
\left[ \sum_i \frac{f_i}{g_i} \right]
$$ 
with $\deg(f_i) = p$ and $\deg(g_i) = q$.
Thus we have
\[
\overline{R^p_q} \subseteq \sum_{X \in L} \overline{R^p_{q, X}}.
\]

Suppose that $\sum_{X \in L} \phi_X =0$ in $\overline{R^p_q}$ 
where $\phi_X \in \overline{R^p_{q, X}}$.
We will prove $\phi_X = 0$ for each $X$.
Fix $X\in L$.  Choose a $\bK$-basis $h_{1},  h_{2}, \dots, h_{m} $ 
for $C_{q, X} $.
Write
\[
\phi_{X} = 
\tau^{p}_{q, X}\left(\sum_{i} [f_{i} ]\otimes h_{i}\right)
=
\left[\sum_{i} f_{i}\, h_{i} \right]
\]
for some $f_{i}\in S^{p}$ for each $i$. 
Since
\[
\left[\sum_{i} f_{i}\, h_{i} \right] =
\phi_{X} 
=
-
\sum_{{Y\in L}\atop{Y\neq X}  } \phi_{Y}, 
\]
we have 
\[
\sum_{i} f_{i}\, h_{i}
\in
R^{p-1}_{q-1}
+\sum_{{Y\in L}\atop{Y\neq X}} S^{p}\, C_{q, Y}.
\]
By Lemma~\ref{3.2}, one has  
$f_{i}\in I(X) $ for each $i$. 
Thus
\[
\phi_X 
= 
\tau^{p}_{q, X}\left(\sum_{i} [f_{i} ]\otimes h_{i}\right)
=
0 .
\ \ \ 
\square 
\]

\bigskip

\section{Proofs}
In this section we prove Theorem \ref{maintheorem},
Corollaries \ref{4.2} and \ref{reflectionarrangement}.
First recall the M\"obius function
$\mu : L \rightarrow \mathbb Z$ in Definition \ref{moebius}.  

\begin{definition}
\label{poincare}
The {\bf Poincar\'e polynomial} of $\A(\Delta)$ 
is defined by
\[
\Poin(\A(\Delta), t) = \sum_{X\in L} \mu(X) (-t)^{\codim X}. 
\]
\end{definition}
The {Poincar\'e polynomial} 
is combinatorially defined
and is known  to be equal to the Poincar\'e polynomial
of 
the topological space
$M(\mathcal{A}(\Delta)) := V \setminus \bigcup_{H\in \A(\Delta)} H$ 
when $\mathbf K = \mathbb C$
\cite{ors1} \cite[Theorem 5.93]{ort1}.

\bigskip
\noindent
{\it Proof of Theorem \ref{maintheorem}.}\/
By Propositions~\ref{3.4} and \ref{conj_main}, we have
\[
\overline{R^p_q} =  \bigoplus_{X \in L} \overline{R^p_{q, X}}
\ \ \ {\rm and} \ \ \
\overline{R^p_{q, X}} \cong  S^p_X \otimes C_{q, X}.
\]
Combining these with
Proposition~\ref{3.1A}, we have
\begin{eqnarray*}
\Poin(\overline{R}(\Delta),s,t) & = & \sum_{p, q} 
\dim(\overline{R^p_q}) s^p t^q 
 =  \sum_{p, q} \left( \sum_{X \in L} 
 \dim(\overline{R^p_{q, X}}) \right) s^p t^q \\
 & = & \sum_{X \in L} \left( \sum_{p \ge 0} \dim(S_X^p) s^p \right)
 \left( \sum_{q \ge 0} \dim(C_{q, X}) t^q \right)\\
 & = & \sum_{X \in L} \left(\frac{1}{1-s}\right)^{\dim X} 
  (-1)^{\codim X} \mu(X) \left( \frac{t}{1-t} 
 \right)^{\codim X} \\
 & = & \sum_{X \in L} (-1)^{\codim X} \mu(X)  
 \left( \frac{t(1-s)}{1-t} \right)^{\codim X} 
 \left(\frac{1}{1-s}\right)^{\ell} \\
 & = & \frac{1}{(1-s)^\ell} \Poin \left(\A(\Delta),\frac{t(1-s)}{1-t} 
 \right). \ \ \ \square
\end{eqnarray*}

\bigskip
\noindent
{\it Proof of Corollary \ref{1.2B}.}
By Propostion \ref{2.2}, we obtain
\[
\dim(R^{p}_{q}) = \sum_{a \le p, b\le q} \dim(\overline{R^{a}_{b}}). 
\]
Thus 
\begin{multline*} 
\sum_{p, q} \dim (R^{p}_{q}) s^{p}t^{q} 
=
\left(
\sum_{c\ge 0} s^{c} 
\right) 
\left(
\sum_{d\ge 0} t^{d} 
\right) 
\sum_{a, b} \dim (\overline{R^{a}_{b}}) s^{a}t^{b} \\
=
\frac{1}{(1-s)(1-t)} \Poin \left(\overline{R}(\Delta), s, t\right)
=
 \frac{1}{(1-s)^{\ell+1} (1-t)} \Poin \left(\A(\Delta),\frac{t(1-s)}{1-t} 
 \right). \ \ \ 
 \square
\end{multline*}

\bigskip

Let $\Der$ be the $S$-module of derivations :
$$
\Der =\{\theta \mid \theta:S \to S \text{ is a } \mathbf{K}\text
{-linear derivations}\}.
$$
Then $\Der$ is naturally isomorphic to $S {\bigotimes}_{\mathbf K}V$.
Define
$$
D(\Delta)=\{\theta \in \Der \mid \theta(\alpha) \in \alpha S 
\text{ for any }\alpha \in \Delta\},
$$
which is naturally an $S$-submodule of $\Der$.
We say that the arrangement $\mathcal{A}(\Delta)$ is {\bf free} if
$D(\Delta)$ is a free $S$-module \cite[Definition 4.15]{ort1}.
An element $\theta\in D(\Delta)$ is said to be {\bf homogeneous of degree
$p$} if $$
 \theta(x)\in S^{p} \text{~for all~}    x \in V^*.
$$
When $\mathcal{A}(\Delta)$ is a free arrangement, let 
$\theta_1, \cdots, \theta_{\ell}$ be a homogeneous basis for $D(\Delta)$.
The $\ell$ nonnegative integers 
$\deg \theta_1, \cdots, \deg \theta_{\ell}$ are called the
{\bf exponents} of $\mathcal{A}(\Delta)$.
Then one has 
\begin{proposition}
{\bf (Factorization Theorem \cite{ter1}, \cite[Theorem 4.137]
{ort1})}
\label{factorizationtheorem}
If $\mathcal{A}(\Delta)$ is a free arrangement with exponents
$d_1, \cdots, d_{\ell},$ then
$$
\Poin(\A(\Delta), t)=\prod^{\ell}_{i=1}(1+d_i t).
$$
\end{proposition}

\noindent
By Theorem~\ref{maintheorem} and Proposition~\ref{factorizationtheorem}, 
we immediately have Corollary~\ref{4.2}. Corollary~\ref{reflectionarrangement}
is a special case of Corollary~\ref{4.2} because the set $\A(\Delta)$ 
of reflecting hyperplanes is known to be a free  arrangement
\cite{sai1} \cite{ter2} \cite[Theorem 6.60]{ort1}. 

\bibliographystyle{unsrt}

\end{document}